\documentclass[preprint,12pt]{elsarticle}
\usepackage[utf8]{inputenc}
\usepackage{graphicx}
\usepackage{amsmath}
\usepackage{tabularx}
\usepackage{array}
\usepackage{float}
\usepackage{makecell}
\usepackage{multirow}
\newcolumntype{M}[1]{>{\centering\arraybackslash}m{#1}}
\newcolumntype{Y}{>{\centering\arraybackslash}X}
\newcolumntype{Z}{>{\raggedright\arraybackslash}X}
\usepackage{amssymb}
\usepackage{tabularx}

\journal{Electric Power Systems Research Journal}

\begin{document}

\begin{frontmatter}

\title{Synergizing Machine Learning with ACOPF: A Comprehensive Overview}

\author[inst1]{Meng Zhao}

\affiliation[inst1]{organization={Energy and Environment Directorate},
            addressline={Pacific Northwest National Laboratory}, 
            city={Richland},
            postcode={PO Box 999, 99352}, 
            state={WA},
            country={USA}}

\author[inst2]{Masoud Barati}

\affiliation[inst2]{organization={Electrical and Computer Engineering, University of Pittsburgh},
            addressline={Swanson School of Engineering}, 
            city={Pittsburgh},
            postcode={15206}, 
            state={PA},
            country={USA}}

\begin{abstract}
Alternative current optimal power flow (ACOPF) problems have been studied for over fifty years, and yet the development of an optimal algorithm to solve them remains a hot and challenging topic for researchers because of their nonlinear and nonconvex nature. A number of methods based on linearization and convexification have been proposed to solve to ACOPF problems, which result in near-optimal or local solutions, not optimal solutions. Nowadays, with the prevalence of machine learning, some researchers have begun to utilize this technology to solve ACOPF problems using the historical data generated by the grid operators. The present paper reviews the research on solving ACOPF problems using machine learning and neural networks and proposes future studies. This body of research is at the beginning of this area, and further exploration can be undertaken into the possibilities of solving ACOPF problems using machine learning.
\end{abstract}



\begin{keyword}
AC optimal power flow \sep machine learning \sep neural network \sep physics-informed \\
\small \textit{This work was supported by the NSF ECCS Award 1711921.}
\end{keyword}

\end{frontmatter}


\section{Introduction}
Optimal Power Flow (OPF), first introduced by Carpentier in 1962 \cite{first}, aims to achieve an optimal operating point in terms of a specified objective function such as minimizing generation cost, minimizing total ohmic losses, matching a desired voltage profile for different generators within a transmission network subject to the physical laws, and operational or technical constraints \cite{Frank2016AnIT}. 
The OPF typically runs over time horizons ranging from a few milliseconds in real time, to 24 hours in advance for power system operation and control for all Independent System Operators (ISO) around the world \cite{Zhao2021ASB}.
The ACOPF is an optimization problem within OPF that considers the full AC power flow equations and can be used for a variety of purposes. The nonlinear power flow equations of ACOPF make these problems nonconvex and NP-hard \cite{r1}. In light of the challenges posed by the nonconvexity and nonlinearity of ACOPF problems, considerable research has been conducted, which can be distinguished by two traditional methods linearization and convexification.

A linearized version of ACOPF is known as direct current optimal power flow (DCOPF), which approximates power flows in a linearized manner while ignoring transmission losses. A comparison between DCOPF and ACOPF is presented in Table \ref{t1} \cite{history}. The advantages of DCOPF include its simplicity, which allows it to solve problems quickly and be applied in large-scale networks. However, there are several challenges associated with DCOPF. Grid conditions may differ from the linear assumptions imposed by DCOPF, so grid failure and instability are more likely to occur \cite {r3}. Relying on DCOPF can also have serious consequences for climate change. As estimated by the Federal Energy Regulatory Commission (FERC) in their 2012 report, approximate-solution techniques can lead to costs of billions of dollars as well as unnecessary emissions \cite{r4}. 


\begin{table}[H]
  \caption{Comparison of DCOPF and ACOPF \cite{history}.}
  \label{t1}
  \scriptsize 
  \begin{tabularx}{\textwidth}{|M{1.5cm}|Y|Y|Y|Y|M{2cm}|Y|Y|}
    \hline
    Problem name & Includes voltage angle constraints? & Includes bus voltage magnitude constraints? & Includes transmission constraints? & Includes losses? & Assumptions & Include generator costs? & Includes contingency constraints? \\
    \hline
    DCOPF & No & No; all voltage magnitudes fixed & Yes & Maybe & Voltage magnitudes are constant & Yes & No  \\
    \hline
    ACOPF & Yes & Yes & Yes & Yes & & Yes & No \\
    \hline
  \end{tabularx}
\end{table}


As a consequence of the inaccuracy and insecurity of DCOPF, it is more practical to work on ACOPF problems directly. Multiple algorithms based on convex optimization have been proposed. The second-order cone programming (SOCP) \cite{Jabr2007ACQ} and semidefinite programming (SDP) \cite{Lavaei2012ZeroDG} have been widely employed and investigated. These studies led to the development of improved algorithms, such as convex iteration and \cite{conv}, which are more efficient. However, these algorithms may not ensure the global optimum \cite{r2} nor guarantee global optimization under precise assumptions, and they may not be compatible with large-scale power systems. Therefore, research on global solutions to ACOPF problems is still ongoing.


As artificial intelligence (AI) and machine learning (ML) have evolved and matured over the past few years, researchers from diverse fields have begun to explore ways to solve complex, unstructured problems by incorporating these new techniques. 

For OPF problems, ML approaches have been proposed to learn the initialization of the OPF problem. An ML-based method was proposed for initializing state estimation within a distribution system \cite{20}. To accelerate the ACOPF solution procedure, \cite{Baker} provides a good starting point for initializing ACOPF. 
Other ML approaches directly bypass typical optimization solvers or iterative methods to learn solutions online. A good example is  \cite{24} which makes use of ML to solve mixed-integer quadratic equations. A feasible and optimal DCOPF solution learning is presented in \cite{25}. 
To predict the active constraints in the ACOPF problem, a ML technique has been developed in \cite{22,23}.  

By using historic or simulation-derived data, \cite{20} trained a shallow neural network to learn the initialization used in Gauss-Newton for distribution system state estimation (DSSE) problems. This combination of ML and traditional optimization approaches shows a superior performance in system stability, accuracy, and computation efficiency. In \cite{23}, neural networks were used to learn the mapping from uncertainty realizations to the active set of a DCOPF problem as an intermediate step towards learning the optimal solution. Once the active set is determined, the optimal solution to the original problem can be recovered by solving a linear system of equations. This work illustrated the benefits of using historic data to exploit the structure of the OPF problem, and solve DCOPF on timescales appropriate for corrective control. In \cite{25}, the authors trained a deep neural network (DNN) model to learn the mapping between the load inputs, the dispatch and transmission decisions for DCOPF problems. They leveraged a useful structure in DCOPF to significantly reduce the mapping dimension, subsequently cutting down the size of the DNN model and the amount of training data/time needed. In addition, a postprocessing procedure was designed to ensure the feasibility of the obtained solutions.

In this body of research, neural networks (NN) demonstrated the ability to model extremely complicated nonconvex functions, making them highly attractive for this setting. A model could be trained offline on historic data and used in real time to make predictions on an optimal power setting \cite{Guha}. In addition, research on solving ACOPF problems with ML algorithms have been widely conducted, and in this paper we will review state-of-the-art methods proposed.

The rest of this paper is organized as follows: Section \uppercase\expandafter{\romannumeral2} briefly introduces the formulation of ACOPF and ML algorithm. Current research about ACOPF with ML is  discussed in Section \uppercase\expandafter{\romannumeral3} and future research directions are discussed in Section \uppercase\expandafter{\romannumeral4}. Section \uppercase\expandafter{\romannumeral5} summarizes this review paper.


\section{Background}

In this section, a general formulation of ACOPF problems, ML algorithms, and artificial neural networks (ANN) are introduced.

\subsection{ACOPF Problem Formulation}

One general ACOPF problem is to minimize energy generation cost while satisfying physical laws and engineering constraints. In this case, the classical problem formulation is as follows:
\begin{equation}{\label{acopf}}
\begin{split}
  & \min \quad \sum_{k \in \mathcal{G}}f_k(P_{Gk})\\
  & \ s.t.\quad \  P_{G_k}-P_{D_k}=V_{dk} \sum_{i=1}^n (G_{ik}V_{di}-B_{ik}V_{qi})+V_{qk} \sum_{i=1}^n (B_{ik}V_{di}+G_{ik}V_{qi}),\\
  & \qquad \ \ \   Q_{G_k}-Q_{D_k}=V_{dk} \sum_{i=1}^n (-B_{ik}V_{di}-G_{ik}V_{qi})+V_{qk} \sum_{i=1}^n (G_{ik}V_{di}-B_{ik}V_{qi}),\\
  & \qquad \ \ \ P_k^{min} \leq P_{G_k} \leq P_k^{max},\\
  & \qquad \ \ \ Q_k^{min} \leq Q_{G_k} \leq Q_k^{max},\\
  & \qquad \ \ \ (V_k^{min})^2 \leq V_{dk}^2+V_{qk}^2 \leq (V_k^{max})^2, k \in \mathcal{N}.
\end{split}
\end{equation}
Where, 

$P_{Gk}$, $Q_{Gk}$ are the active and reactive power of generators, respectively;

$P_{Dk}$, $Q_{Dk}$ are the active and reactive power of demands, respectively;

$V_{dk}$ and $V_{qk}$ are the real and imaginary parts of voltage $V_k$, respectively;

$G_{ik}$ and $B_{ik}$ are the real part and imaginary part of the admittance matrix $Y$, respectively.

Given the active and reactive power demands at the load buses, the active and reactive power generated by the generators and voltage magnitudes at each bus can be obtained by solving \eqref{acopf}. From \eqref{acopf}, it is clear that ACOPF problems are nonlinear and nonconvex, indicating the difficulty and impossibility of obtaining  global optimal solutions directly.
Over the last decade, numerous attempts have been made by researchers in academia to solve the nonconvex OPF problem through a variety of convex relaxation methods, such as second order conic (SOC) relaxation \cite{6756976}, semidefinite programming (SDP) \cite{5971792}, and quadratic convex relaxation \cite{7271127}. Although these methods generated initial excitement, they failed to properly handle large-scale systems. One particular drawback of these methods is the difficulty in recovering a feasible solution \cite{6980142}. Another research trend has investigated reformulation of the problem using rectangular current and voltage, instead of the common formulation with power and polar voltage \cite{onil}. While the network current flow in this formulation becomes linear, some of the other linear terms will become nonlinear and nonconvex \cite{7299702}. The nonlinear terms and the nonconvexities of the feasible space can potentially be handled using a successive linear programming approach \cite{7299702}. However, due to the nature of the formulation, decoupling active and reactive power is not possible. 
A few disadvantages are associated with the methods discussed in this article: They rely heavily on mathematical optimization and do not take advantage of the characteristics of electric power flow. And despite recent advancements in computing, these methods are not yet computationally tractable and cannot handle large-scale systems within time frames available.

\subsection{Machine Learning Algorithm}

\begin{figure}[H] 
\centering
\centerline{\includegraphics[width=0.3\textwidth]{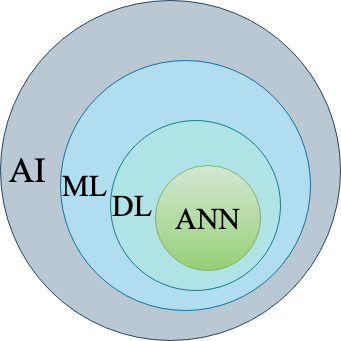}}
\caption{Relationship among AI, ML, DL and ANN.}
\label{rela}
\end{figure}

A subset of AI, ML is the study of computer algorithms that improve automatically through experience \cite{ml}. To perform ML, the first step is to train a model with training data and then use the model to make predictions on other similar data. Within ML, deep learning (DL) is the subfield of ML and ANN (or more simply neural network (NN)) is the subfield of DL. Fig. \ref{rela} shows the relationship among them.

ML algorithms utilize data mining to inform future models with historic data, without knowing the causal relationships or physical laws. Basically, there are three steps in ML algorithms: Step 1) the decision process, in which a ``guess" is made based on the pattern of the data; Step 2) an error function, a measuring method that judges the quality of the “guess” by comparing it with known results; Step 3) the optimization process, where the algorithm aims to minimize the error in Step 2 and improves the decision process until it becomes the final decision \cite{mll}. Based on its intelligence and versatility, ML techniques have been used in a wide variety of fields, including image recognition, autonomous driving, the energy sector, and health care. Furthermore, a growing number of tech companies are utilizing algorithms like these to provide, for instance, their clients with more relevant advertisements (Amazon) or to estimate delivery times (Uber Eats) to provide a better user experience. ML can now solve many previously unsolvable problems. 


\subsection{Artificial Neural Network (ANN)}

\begin{figure}[H] 
\centering
\centerline{\includegraphics[width=0.6\textwidth]{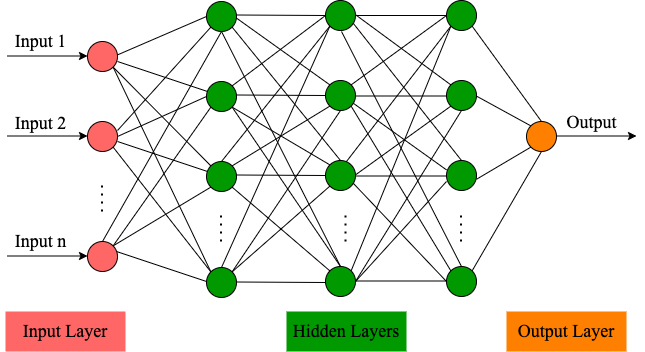}}
\caption{Artificial Neural Network.}
\label{ann}
\end{figure}

 Inspired by the biological neural networks in animal brains, ANN has been widely used. An NN model is composed of the connected units or nodes called ``neurons". The neurons can receive and process a signal (data) and then tranmit it to the following neuron through the connection between them, called the ``edges." During the training process, the weights (which denote the strength of each signal (data) in the neurons) and edges will be updated. A general concept of ANN is shown in Fig. \ref{ann}.

A unique characteristic of ANN is that they can learn the complex relationships between dependent and independent variables without assumptions about their mathematical representation. Therefore, the ANN is a good tool to model complicated (non)linear relationships and nonconvexity between input and output datasets, and it has been used across a range of disciplinary applications (e.g., detecting neutrinos \cite{1}, climate simulations \cite{2}, and fluid dynamic simulation \cite{3}, etc.). Since ANN can reproduce and model nonlinear processes, it has been used in various areas: quantum chemistry \cite{che}, general game playing \cite{game}, signal classification \cite{sc}, lung cancer diagnosing \cite{cancer}, and etc. It is noticeable that ANN has a high flexibility and can be used in multiple disciplines. With ANN, researchers can obtain more accurate results with reduced latency \cite{4}.


\section{ACOPF with Neural Network}
Given the benefits of using ANN to solve diverse problems in other domains, investigating this algorithm to solve ACOPF problems is appropriate. Recent research on ACOPF with NN can be classified into four directions: 1) train a model offline with historic data to learn the mapping between system loads and optimal generation set values  to find a near-optimal and feasible ACOPF solution without actually solving an optimization problem \cite{Zamzam2020LearningOS}; 2) generate NN models to predict the initial values of variables critical to the problem's convergence as the high-quality initial condition for existing power-grid numerical solvers \cite{pgsim}; 3) utilize deep neural network (DNN) to solve ACOPF problems directly, and combine  this strategy other methodologies to improve the computation efficiency; 4) develop physics-informed NN which considers the constraints of the system while training the model to improve its robustness without violating the constraints. Except for ACOPF problems, ML algorithms were also leveraged in their byproduct problems within OPF field.

\subsection{Shallow Neural Network}

Within the power systems realm, grid operators have generated a large quantity of data by repeatedly solving ACOPF problems throughout the day \cite{Zamzam2020LearningOS}. These data can be used by ML algorithms to train and test the NN developed for the ACOPF problems. 

\cite{pred} discussed predicting the output of ACOPF problems by ML algorithms. Based on the data generated from the IEEE 30-bus power system, they compared the training and test results of four algorithms in Python: stochastic gradient descent, support vector regression, multi-layer perceptron regression and gradient boosting regression (GBR). Among these methods, GBR performed the best and computation time was reduced by up to 30\%. However, approximately 60\% of the predictions violated one of the constraints in the ACOPF problems even though the GBR algorithm can predict the power and voltage generated within high precision. Also, they did not test on larger-scale power systems and each output was predicted independently, which ignored the relationship between voltage and power.

\begin{figure}[H] 
\centering
\centerline{\includegraphics[width=0.6\textwidth]{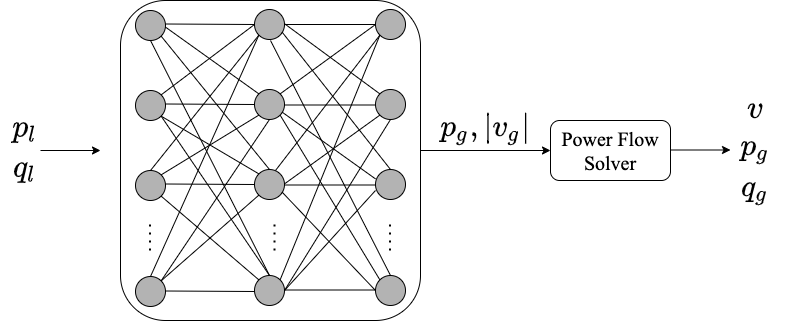}}
\caption{Workflow of the proposed learning approach in \cite{Zamzam2020LearningOS}.}
\label{fastflow}
\end{figure}

Inspired by \cite{pred}, the authors in \cite{Zamzam2020LearningOS} utilized ML to learn the mapping between network loads and optimal power generations. An NN model was trained offline and was then used to solve the ACOPF problems online. To train the NN model, strictly feasible solutions were used and the outputs of the model were constrained to satisfy the active power generation and voltage magnitude constraints. Fig. \ref{fastflow} shows the architecture of the proposed approach. The ultimate framework was faster and provided a more optimal solution than the traditional ACOPF approximation algorithms.

In \cite{Baker2020ALQ}, a quasi-Newton method based on ML was presented. It performed iterative updates for candidate optimal solutions without having to calculate a Jacobian or approximate Jacobian matrix. The proposed learning-based algorithm developed an NN with feedback and it is capable of finding approximate solutions to ACOPF very quickly. With the proper choice of weights and activation functions, the model became a contraction mapping, and convergence of the solutions can be guaranteed.



By training an NN model with historical data, neither linearization nor convexification of the ACOPF problem is needed. Also, this algorithm is less computationally burdensome and is more efficient than traditional methods. However, results may violate the constraints of the problems since they are not considered during the training process.

\subsection{Learning Warm Starts for ACOPF Solvers}

Except training a model to predict the solutions of ACOPF problems mentioned above, some researchers sought to accelerate power grid numerical solvers by using ML algorithms to learn the warm starts of the solvers like MIPS (MATPOWER Interior Point Solver), which can guarantee convergent solutions and improve computational efficiency.

MIPS is a solver in MATPOWER that can solve power flow and OPF problems \cite{MIPS}. First, MIPS converts the inequality constraints in \eqref{acopf} into equality constraints. Second, MIPS uses a Lagrangian formulation to reformulate the ACOPF problem and leverages the Newton method to solve the reformulated Lagrangian problem. However, the Newton method has a heavy computational burden and is time consuming since it iteratively converges to a set of convergence criteria. Therefore, improving the efficiency of this kind of ACOPF solver is necessary and practical.

\begin{figure}[H]
\centering
\centerline{\includegraphics[width=0.6\textwidth]{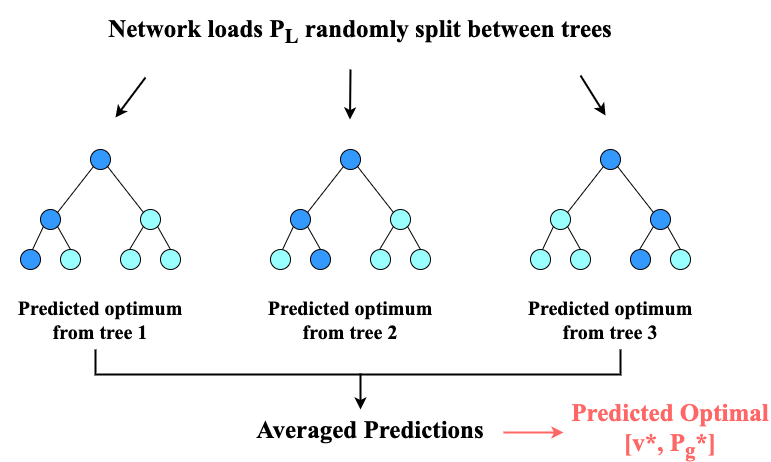}}
\caption{Workflow of the proposed learning approach in \cite{Baker}.}
\label{rf}
\end{figure}

\cite{Baker} put forward a multitarget approach to train a Random Forest model \cite{random} which is shown in Fig. \ref{rf} with historic optimal power flow solutions. 

In this algorithm, only the power demands, the number of generators, and buses in the power system were needed; the network topology did not need to be known (i.e., the admittance matrix, line parameters, etc.). The inputs of the model are the loads at each bus, the outputs are the voltage magnitudes and optimal power generations at each bus. Since the voltages and power generations are correlated with each other, this multi-target model can simultaneously predict multiple outputs. In this model, each Decision Tree agent performs its individual regression procedure separately, then  at the end the agents coordinate these outputs by averaging them. This model can help alleviate the risks of over-fitting associated with individual Decision Trees. The learned warm start using the MIPS solver can get a convergent solution faster than a DC warm start or a flat start.

To speed up computations, \cite{pgsim} proposed a framework named "Smart-PGSim" to estimate a high-quality initial start for the MIPS solver which can guarantee the precision and robustness of the solutions. Unlike prior methods, all of the inputs in AC-OPF problems were considered in the algorithm and the predicted solution was guaranteed optimal while the efficiency and performance improved.

\begin{figure*}[h!]
\centerline{\includegraphics[width=1\textwidth]{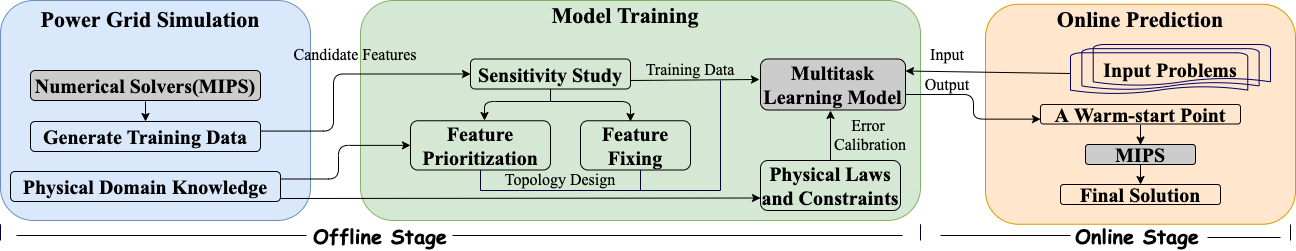}}
\caption{Workflow of Smart-PGSim in \cite{pgsim}.}
\label{spg}
\end{figure*}

The proposed Smart-PGsim framework in \cite{pgsim} is shown in Fig.~\ref{spg} which included two phases: offline and online. In offline phase, it finds the most important features to construct an efficient NN model for online acceleration by using the power grid simulation. That is, during the sensitivity study of the offline phase, it helps identify which variables in MIPS are necessary as the output of the prediction model and quantifies the impact of the imprecise variables (i.e., variables with some accuracy loss) on the success rate of simulation and performance in terms of execution time. These results are then used to help design the multi-task learning (MTL) model. In online phase, the MIPS solver can use the high-quality start points generated by the MTL model for quick convergence.


These proposed algorithms can generate a warm start of the MIPS solver or other conventional solvers of ACOPF, and they can help the solvers to reach a more convergent solution in less time.

\subsection{Deep Neural Network (DNN)}

As a subfield of ML, DL can structure algorithms in layers and create an ANN which can learn and make intelligent decisions on its own \cite{DeepLM}. Instead of considering limited constraints and variables of the ACOPF, or learning warm starts of the ACOPF solvers, some research sought to solve the ACOPF problems directly by using state-of-the-art DL algorithms.

\cite{Zhou2020ADM} derived fast ACOPF solutions with the deep reinforcement learning (DRL) algorithm in which imitation learning was used to generate initial weights for the NN, and a proximal policy optimization (PPO) algorithm \cite{ProximalPO} was leveraged to train and test the AI agents. In addition, the original ACOPF problem was formulated as Markov decision process with appropriate power grid simulation environment and a reasonable reward function, so that the DRL agent can be trained without violating the operation constraints under various operation conditions. The computation time was improved at least by $7$ times compared with conventional interior-point solver, and the model also worked under $N-1$ contingency on changing network topology.

By combining with DNN and Levenberg-Marquardt back-propagation, \cite{RealTimeOP} proposed the twin delayed deep deterministic policy gradient (TD3) approach to improve the computational performance of the ACOPF problems when the demands in the system have a random pattern. Also, random Gaussian noise was added in individual net loads to represent the uncertainty characteristics introduced by renewable energy sources. During the training process, the appropriate reward vector was set based on constraints, and solutions of the final model were compared with the MATPOWER solutions, which implied its robustness and efficiency. 

In \cite{Singh2021LearningTS}, the authors trained a sensitivity-informed DNN (SI-DNN) to match not only the OPF optimizers, but also their partial derivatives with respect to the OPF parameters (loads). The training data was extended with the Jacobian matrix which carries the partial derivative of the minimizer and the Forbenius norm of the matrix was added in the loss function to incorporate sensitivity information into the training process. This sensitivity-informed training can readily complement other existing learn-to-OPF methodologies, and the proposed SI-DNN model offers the same prediction performance with less training data than conventionally trained DNN models, which can improve the sampling efficiency.


The research discussed combined DNN and other algorithms to solve ACOPF problems directly in which the prediction accuracy and computation efficiency were improved significantly. Some trained models are even robust to system typologies and can be combine with existed methodologies.

\subsection{Physics-Informed Neural Network}

The research introduced in the previous subsections developed NN models of ACOPF problems only with generated load demands, without knowing the topology or considering constraints in the problems. This limitation can cause a less-robust NN model and the solutions may validate physical constraints. To tackle this issue, some researchers have worked on physics-informed NN \cite{Physicsin} for the ACOPF problems in which the constraints were added in the loss function, or the violation degree of the constraints were defined during the training process. In this case, the feasibility of the solutions are guaranteed without violating the constraints.

The Smart-PGSim developed in \cite{pgsim} also utilized physical domain knowledge which denotes the physical constraints of the system during the training process. Power balance equations, inequality constraints, and cost function were all considered and added into the loss function. Therefore, the trained model has higher accuracy, feasibility, and robustness.

\begin{figure*}[h!]
\centerline{\includegraphics[width=1.0\textwidth]{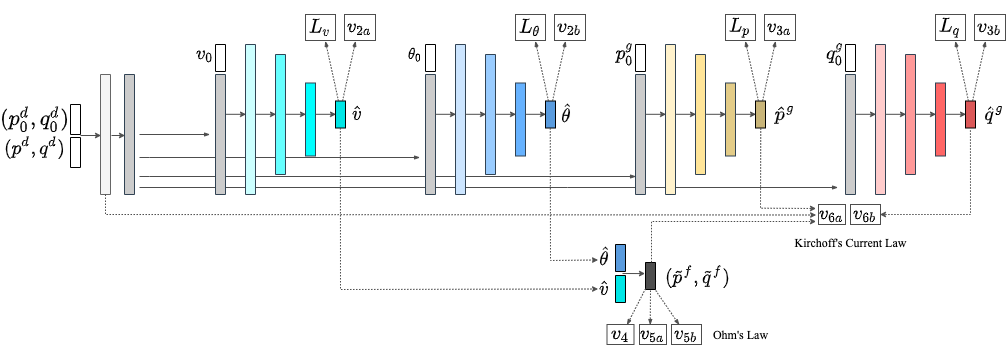}}
\caption{The structure of OPF-DNN model in \cite{dnn}.}
\label{dn}
\end{figure*}

In \cite{dnn}, to predict generator set points, the authors designed a DNN model to the ACOPF (OPF-DNN) which is shown in Fig. \ref{dn}. Also, the model leveraged the physical and engineering constraints in a Lagrangian framework with violation degrees by extending the predicted values with the reactive power dispatch and the voltage angles. In the OPF-DNN model, each layer is fully connected with the ReLU activation and each layer takes the outputs of the previous layer as its inputs. This algorithm enhanced the prediction accuracy and efficiency compared with previous methods.

In \cite{DeepOPFAF}, an algorithm named DeepOPF was proposed to solve ACOPF problems in which a DNN model was trained first to predict the independent operating variables; the remaining dependable variables were then obtained by solving the AC power flow equations. A penalty term with zero-order estimation technique, which captured the violation of the inequality constraints, was added in the loss function during the training process. However, this model is not robust to  changes in system typologies or parameters.

\begin{figure}[H]
\centering
\centerline{\includegraphics[width=0.6\textwidth]{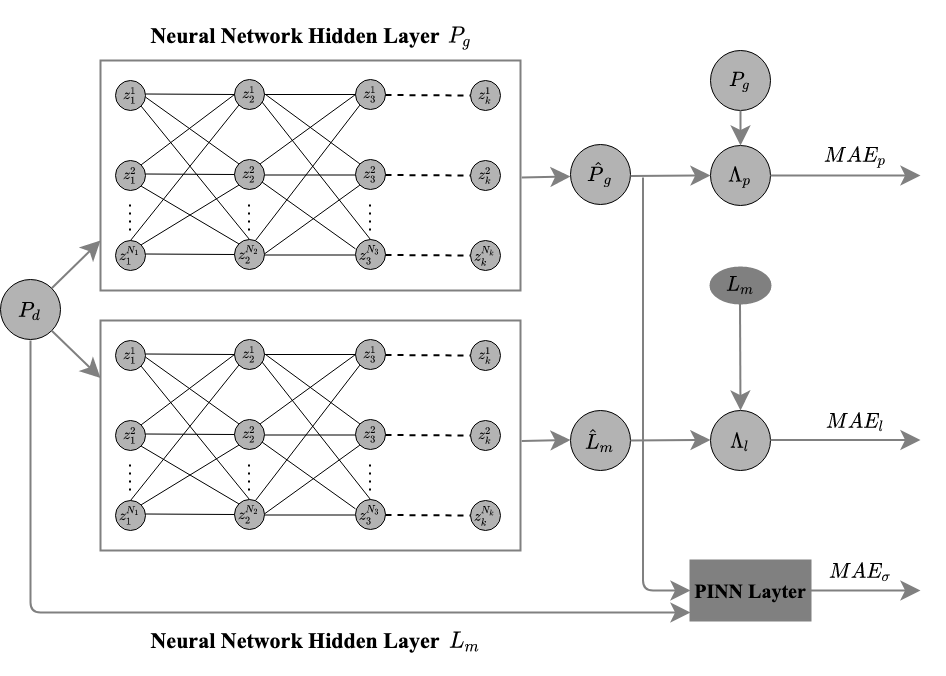}}
\caption{Structure of the physics-informed NN in \cite{PhysicsInformedNN}.}
\label{dcopf}
\end{figure}

To obtain more accurate results with fewer data, \cite{PhysicsInformedNN} designed a physics-informed NN with worst-case guarantees for the DCOPF problems. In this proposed physics-informed NN, shown in Fig. \ref{dcopf}, the physical equations of the DCOPF problem were incorporated into the neural network loss function and the generated optimal value should satisfy the KKT conditions, which were set as necessary and sufficient conditions. Moreover, the authors tried to find the worst-case violations of the physics-informed NN, and developed methods to reduce them. Test results showed that the proposed physics-informed NN can provide a more accurate result with fewer data, and this process can be used to recognize both the potential and the  challenges for its future application to ACOPF problems.

The research introduced in this subsection trained the physics-informed NN, which considered the constraints of the ACOPF problems, resulting in a more robust and accurate model. 

\subsection{Test Results and Evaluation of the ML Models}

To validate the effectiveness of the proposed ML models, the studies presented performed simulations on select power systems. The IEEE 118-bus system was widely used in these papers and Table \ref{result} shows the corresponding results of each ML model. From Table \ref{result}, it is clear that solving ACOPF problems with ML algorithms can significantly improve computational efficiency, and the model trained with more information on the system can generate more accurate predictions.

\begin{table*}[h!]
  \centering
  \caption{Test Results of Current ML Algorithms on IEEE 118-bus system.}
  \label{result}
  \scriptsize
  \begin{tabularx}{\textwidth}{|c|Y|Y|Y|Y|}
    \hline\hline
    \makecell{Research \\ Cited} & \makecell{Speed-up \\ (with definition)} & \makecell{Prediction \\ Error/MSE} & Optimality & Feasibility \\
    \hline
    \cite{Zamzam2020LearningOS} & \makecell[c]{$11.83X$ \\ $\left(\frac{1}{T}\sum_{t=1}^{T}\frac{\tau_{\text{solver},t}}{\tau_{\text{ML},t}}\right)$} & $\times$ & $2.974 \times 10^{-4}$ & \makecell[c]{$1.41 \times 10^{-8}$ \\ (infeasibility)} \\
    \hline
    \cite{Baker} & $\times$ & $2.12\%$ ($p_g$); $0.01\%$ ($v$) & $\times$ & $\times$ \\
    \hline
    \cite{pgsim} & \makecell[c]{$3.28X$ \\ faster than the \\ traditional \\ MIPS solver} & $0\%$ & $\times$ & $\times$ \\
    \hline
    \cite{RealTimeOP} & $\times$ & \makecell[c]{$0\%$ \\ compared with \\ MATPOWER \\ MIPS} & $\times$ & $\times$ \\
    \hline
    \cite{Singh2021LearningTS} & \makecell[c]{$63X$ \\ faster than \\ conventional \\ OPF solvers} & $0.9$ (test MSE) & $\times$ & $\times$ \\
    \hline
    \cite{dnn} & $>10^4X$ & \makecell[c]{$0.034\%$ ($p_g$); \\ $0.618\%$ ($q_g$); \\ $0.029\%$ ($v$); \\ $0.207\%$ ($\theta$); \\ $0.455\%$ ($p_f$)} & $\times$ & $\times$ \\
    \hline
    \cite{DeepOPFAF} & \makecell[c]{$25X$ \\ faster than \\ Pypower \\ \cite{pypower}} & \makecell[c]{$<0.1\%$ \\ (cost difference \\ with conventional \\ solver)} & \makecell[c]{$<0.1\%$ \\ (optimality \\ loss)} & $100\%$ \\
    \hline
  \end{tabularx}
\end{table*}


Besides discovering multiple ML algorithms to solve ACOPF problems, this body of research also adopted specific metrics to evaluate performance: 1) speed-up factor: a fraction between the time consumed by existing conventional ACOPF solvers and the time consumed by the proposed ML models; 2) feasibility: an indix obtained by evaluating the solutions' satisfaction of the equality and inequality constraints of the ACOPF problems; 3) optimality: this metric is the comparison between the obtained solutions of the proposed ML algorithms and the solutions of the current ACOPF solvers. Moreover, prediction evaluation, which denotes the model's predictive capability, is also an important factor of the proposed network, such as mean absolute error (MAE), mean absolute percentage error (MAPE) and prediction error. Table \ref{indices} summarizes these indices presented in papers cited and the performance of the proposed NN models.

\begin{table}[]
\centering
\caption{Summarized Performance of Current ML Algorithms on ACOPF Problems}
\label{indices}
\begin{tabular}{|c|c|ccc|ccc|}
\hline
\multirow{2}{*}{\begin{tabular}[c]{@{}c@{}}NN \\ Category\end{tabular}}                           & \multirow{2}{*}{\begin{tabular}[c]{@{}c@{}}Research \\ Cited\end{tabular}} & \multicolumn{3}{c|}{ACOPF Metrics}                                              & \multicolumn{3}{c|}{Prediction Matrics}                                               \\ \cline{3-8} 
                                                                                                  &                                                                            & \multicolumn{1}{c|}{Speed-up}     & \multicolumn{1}{c|}{Fes.}  & Opt.   & \multicolumn{1}{c|}{MAE}          & \multicolumn{1}{c|}{MAPE}         & Error \\ \hline
\multirow{3}{*}{\begin{tabular}[c]{@{}c@{}} \\ SNN \\ \end{tabular}}             &  \cite{pred}                                               & \multicolumn{1}{c|}{$\times$}     & \multicolumn{1}{c|}{$\times$}     & $\times$     & \multicolumn{1}{c|}{$\times$}     & \multicolumn{1}{c|}{$\times$}     & $\checkmark$     \\ \cline{2-8} 
                                                                                                  &  \cite{Zamzam2020LearningOS}                               & \multicolumn{1}{c|}{$\checkmark$} & \multicolumn{1}{c|}{$\checkmark$} & $\checkmark$ & \multicolumn{1}{c|}{$\times$}     & \multicolumn{1}{c|}{$\times$}     & $\times$         \\ \cline{2-8} 
                                                                                                  &  \cite{Baker2020ALQ}                                       & \multicolumn{1}{c|}{$\times$}     & \multicolumn{1}{c|}{$\checkmark$} & $\checkmark$ & \multicolumn{1}{c|}{$\checkmark$} & \multicolumn{1}{c|}{$\checkmark$} & $\times$         \\ \hline
\multirow{2}{*}{\begin{tabular}[c]{@{}c@{}}LWS\\ ACOPF \end{tabular}} &  \cite{Baker}                                              & \multicolumn{1}{c|}{$\times$}     & \multicolumn{1}{c|}{$\times$}     & $\times$     & \multicolumn{1}{c|}{$\times$}     & \multicolumn{1}{c|}{$\times$}     & $\checkmark$     \\ \cline{2-8} 
                                                                                                  &  \cite{pgsim}                                              & \multicolumn{1}{c|}{$\checkmark$} & \multicolumn{1}{c|}{$\times$}     & $\times$     & \multicolumn{1}{c|}{$\times$}     & \multicolumn{1}{c|}{$\times$}     & $\checkmark$     \\ \hline
\multirow{3}{*}{\begin{tabular}[c]{@{}c@{}} \\ DNN \\  \end{tabular}}       &  \cite{Zhou2020ADM}                                        & \multicolumn{1}{c|}{$\checkmark$} & \multicolumn{1}{c|}{$\checkmark$} & $\times$     & \multicolumn{1}{c|}{$\times$}     & \multicolumn{1}{c|}{$\times$}     & $\times$         \\ \cline{2-8} 
                                                                                                  &  \cite{RealTimeOP}                                         & \multicolumn{1}{c|}{$\times$}     & \multicolumn{1}{c|}{$\times$}     & $\times$     & \multicolumn{1}{c|}{$\checkmark$} & \multicolumn{1}{c|}{$\times$}     & $\times$         \\ \cline{2-8} 
                                                                                                  &  \cite{Singh2021LearningTS}                                & \multicolumn{1}{c|}{$\times$}     & \multicolumn{1}{c|}{$\times$}     & $\checkmark$ & \multicolumn{1}{c|}{$\times$}     & \multicolumn{1}{c|}{$\times$}     & $\times$         \\ \hline
\multirow{3}{*}{\begin{tabular}[c]{@{}c@{}} \\ PINN \\ \end{tabular}} &  \cite{dnn}                                                & \multicolumn{1}{c|}{$\checkmark$} & \multicolumn{1}{c|}{$\checkmark$} & $\times$     & \multicolumn{1}{c|}{$\times$}     & \multicolumn{1}{c|}{$\times$}     & $\times$         \\ \cline{2-8} 
                                                                                                  & \cite{DeepOPFAF}                                              & \multicolumn{1}{c|}{$\times$}     & \multicolumn{1}{c|}{$\checkmark$} & $\times$     & \multicolumn{1}{c|}{$\checkmark$} & \multicolumn{1}{c|}{$\checkmark$} & $\times$         \\ \cline{2-8} 
                                                                                                  & \cite{PhysicsInformedNN}                                      & \multicolumn{1}{c|}{$\times$}     & \multicolumn{1}{c|}{$\checkmark$} & $\checkmark$ & \multicolumn{1}{c|}{$\times$}     & \multicolumn{1}{c|}{$\times$}     & $\times$         \\ \hline
\end{tabular}
\tiny{SNN: Sallow Neural Network; LWS ACOPF: Learning Warm Start for ACOPF Solver; DNN: Deep Neural Network; PINN: Physics-Informed Neural Network; Speed: Speed-up; Fes.: Feasibility; Opt.: Optimality; Error: Prediction Error} \\
\end{table}

\subsection{Byproducts of ACOPF with Machine Learning}

Except for ACOPF problems, ML algorithms were also utilized in other problems within OPF filed by some researchers, such as state estimation, load forecasting, optimal load shedding, electricity market and so on.

\subsubsection{Load Forecasting}

\begin{figure}[H]
\centering
\centerline{\includegraphics[width=0.6\textwidth]{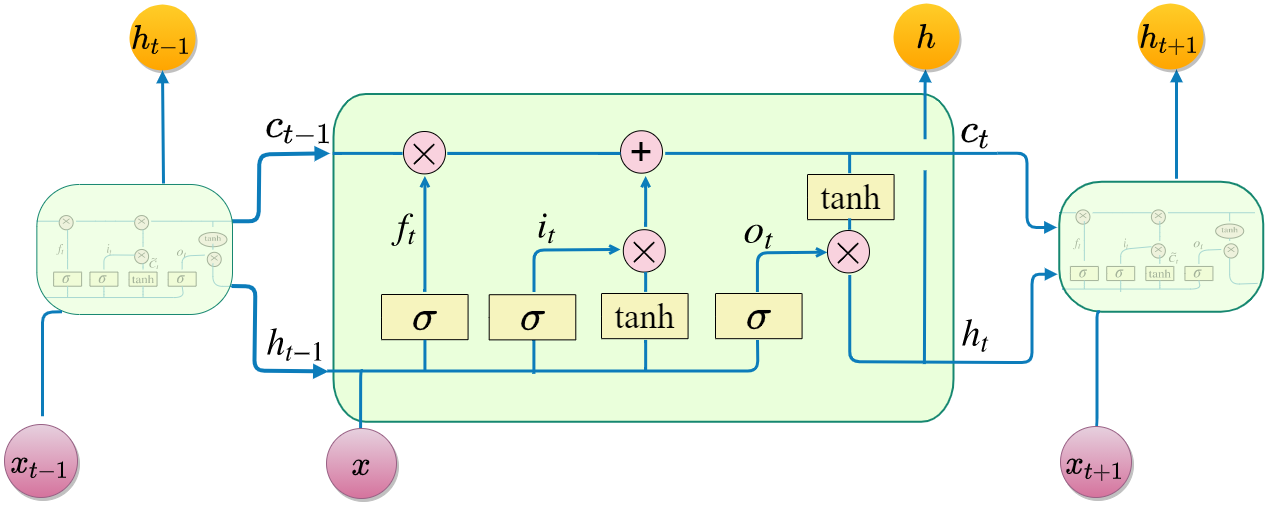}}
\caption{Diagram of the LSTMs.}
\label{lstm}
\end{figure}

In \cite{2019ShortTermRL}, a long short-term memory (LSTM)-based framework was proposed for the short-term load forecasting of individual electric customers in the power system with higher penetration of renewable energy generation. Compared with recurrent neural network (RNN), LSTM can overcome the gradient vanishing issue with a memory cell and forget gate. Fig. \ref{lstm} shows the structure of the LSTM algorithm. In addition, the differences between load on a system level and loads in individual households were compared, and the inconsistency in the residential load profiles was also evaluated to demonstrate the difficulty of residential load forecasting. The inputs of the LSTM model are four vectors which represent the residential energy consumption and time/day indices, respectively. The outputs of the LSTM model are the predicted energy consumption at a specified time interval. Compared with other conventional methods, the LSTM-based model performed best for residential load forecasting.

\begin{figure}[H]
\centering
\centerline{\includegraphics[width=0.6\textwidth]{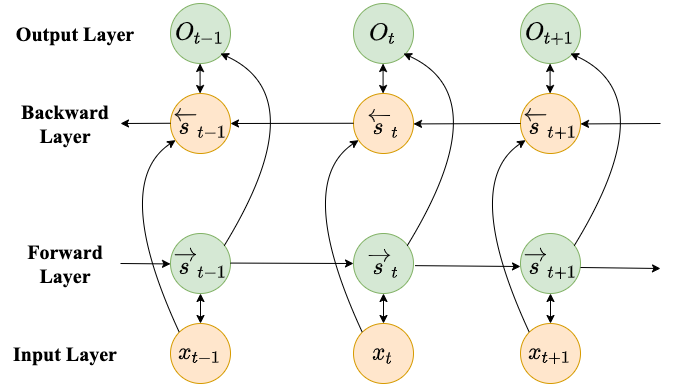}}
\caption{Diagram of the bidirectional LSTMs in \cite{bilstm}.}
\label{bilstm}
\end{figure}

To improve the forecasting capability of the time-series problem, \cite{bilstm} proposed a multi-layer stacked bidirectional LSTM model for load forecasting in power systems. By integrating reverse computing with forward computing, the bidirectional LSTM framework can avoid data loss in long-distance transmission and improve the dependency and reliability of the data. A brief structure of the bidirectional LSTM is shown in Fig. \ref{bilstm}. Test results validated the effectiveness of the proposed model: Higher performance and the capability to retain more original information are guaranteed when compared with BP neural network, ELM, and traditional LSTM. 

To increase the accuracy of short-term load forecasting, \cite{cnnlstm} combined convolutional neural network (CNN) with LSTM;  the model was adopted in Bangladesh power system to predict short-term electrical load. Compared with other existing techniques such as LSTM and extreme gradient boosting algorithm, the proposed framework had higher precision and accuracy.

\subsubsection{Load Shedding}
\cite{Kim2019GraphCN} used graph convolutional network (GCN) to predict an optimal load-shedding ratio that prevents transmission lines from being overloaded under line contingency. In addition, the topology of the power system was convoluted over the NN. The proposed GCN model outperformed both the classical NN model and a linear regression model. 

\subsubsection{State Estimation}

State estimation of unobservable distribution systems was considered in \cite{BayesianSE} and a Bayesian state estimation algorithm based on DL was proposed. They adopted a generic static power flow model of a three-phase power system and used distribution learning for stochastic power injection. During the DNN training process, a Monte Carlo technique was used to compute the minimum mean-squared-error (MMSE) estimate directly from the measurements. Also, a Bayesian bad-data detection and filtering algorithm was developed to detect and remove bad data before the state estimation. Simulation results illustrated the robustness of the proposed algorithm, which worked better than pseudo-measurement techniques on existing benchmarks.

\subsubsection{Unit Commitment}

\cite{2021ModelingTA} replaced the AC power flow equations with a NN-based piece-wise linear model which was then transformed into mixed integer linear programs (MIPS) and embedded into the unit commitment (UC) problem. Moreover, iterative bound tightening, ReLU pruning, and parameter matrix sparsification were utilized to compress the NNs of the model, resulting in the reduced number of binary variables in the MIPS. Test results showed the feasibility and effectiveness of the proposed NN-based model when compared with the solutions of DC and linearized power flow approximations. 

\subsubsection{Electricity Market}
In \cite{GraphNN}, the authors proposed a topology-aware graph neural network (GNN) to predict the electricity market prices by solving OPF problems. The connection between grid topology and locational marginal prices (LMPs), as the outputs of the OPF problem were also studied to improve topology adaptivity and reduce the model complexity. 

\subsubsection{Reliability Assessment}

Probabilistic reliability assessment in operation planning plays an essential role in guiding the operation planner during the decision-making period \cite{2018UsingML}. In \cite{2018UsingML}, the authors combined Monte Carlo simulation, ML, and variance reduction techniques to accelerate the crude Monte Carlo method. The ML algorithm was used to replace the security-constrained OPF (SCOPF) computations with a faster proxy of real-time operator response. The three-area IEEE-RTS96 benchmark was tested by the proposed approach; it improved computation efficiency without sacrificing accuracy.

Moreover, \cite{2020RecentDI} reviewed recent research on the reliability assessment and control problems with ML techniques in power systems. Increasingly, studies are leveraging ML methods to solve these problems, and the proposed algorithms have greate protential to improve the computation process of the reliability assessment. 

\subsubsection{Other Related Research}

\cite{InverseOP} utilized a neural network to formulate the inverse optimal power flow (inverse OPF) algorithm, which can learn private quantities from public quantities. In this paper, the authors tried to recover both cost and admittance matrix parameters by solving the inverse OPF problem via ACOPF with publicly available data, such as five-minute electricity prices and hourly power outputs of electricity generators. From this research, it is noticeable that which critical power grid information is exposed by the published data, also the private and public data are related via an optimization problem. 

To solve probabilistic optimal power flow (POPF) under renewable and load uncertainties of arbitrary distribution, \cite{GaussianPL} proposed a novel Gaussian process learning-based probabilistic optimal power flow (GP-POPF) algorithm which utilizes a nonparametric Bayesian inference-based uncertainty propagation approach. This method needs fewer numbers of samples and less computation time; simulation results showed the generality of the model--- it can work with different types of input uncertainty distributions.

The literature listed in this subsection utilized ML methods to solve OPF-related problems in power systems, which illustrate the efficiency and universality of ML algorithms. 

\section{Future Research Directions}

State-of-the-art research on ACOPF and its byproducts using ML have been discussed. From this body of research, it is apparent that AI approaches including machine (deep) learning are promising and appropriate tools to solve complicated optimization problems. However, this work is still at an early stage of exploration, and more research is needed to make these efforts  applicable to power system industries. Inspired by this current research, future research directions on ACOPF with ML algorithms are discussed.

\subsection{Solving the Relaxed ACOPF with ML Algorithms}

From Section \uppercase\expandafter{\romannumeral3}, it is obvious that current research mainly focuses on solving the ACOPF problems directly by using ML or DL, and most efforts are based on learning the topology of the ACOPF problems; this may lack of sensitivity when some parameters are changed in the system. 

Except considering the original ACOPF problems, one possible line of future research may focus on relaxed ACOPF problems, such as the SDP-relaxed format which is to get the rank-1 solutions. Suppose the $vv^T$ items in the constraints of \eqref{acopf} are relaxed by $W=vv^T$, the final problem is to obtain the rank-1 solutions of the semidefinite symmetric matrix $W$. To solve this relaxed problem, some current research has used conventional algorithms \cite{Zhao2021ASB,Zhao2021LoworderMR}, but none of these strategies recover the exact rank-1 solutions of $W$. Therefore,  ML algorithms might be leveraged to recover the rank-1 $W$ matrix by learning the correlations of each item in $W$ as well as the constraints of the system. There are also some algorithms based on rank minimization problems, such as the penalty method in \cite{Shen2018APM}. Another possibility is to adopt ML algorithms in these existing methods to improve their efficiency and accuracy. 

\subsection{Improving Current ACOPF Solvers with ML Algorithms}

As discussed in Part $B$ of Section \uppercase\expandafter{\romannumeral3}, some current research is leveraging ML algorithms to generate warm starts for existing ACOPF solvers like MIPS to save computation time. Although these methods have improved computational efficiency to some extent, there are still drawbacks: 1) iterations are essential in the solvers to get the final results; and 2) there is no guarantee that whether the warm starts are good or bad (good warm starts can lead to a feasible convergent solution, but bad starts will result in an infeasible solution which is useless and waste of time). 

\begin{figure}[H]
\centering
\centerline{\includegraphics[width=0.6\textwidth]{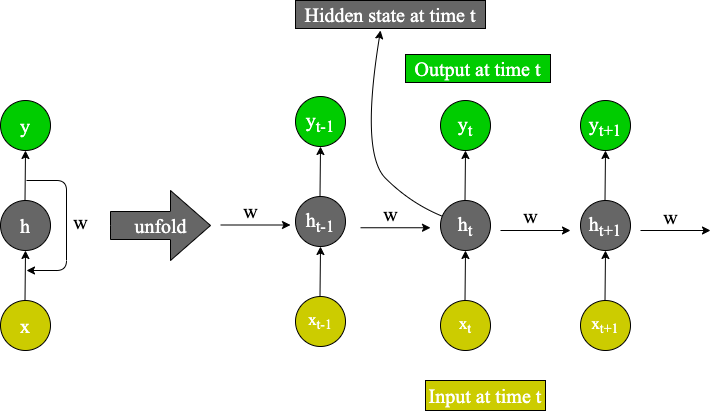}}
\caption{Diagram of the RNN algorithm.}
\label{rnn}
\end{figure}

Since ML algorithms are good at learning and adapting, the iteration process in ACOPF solvers can be learned and substituted by them. For example, RNN in which the outputs depend on the prior elements within the sequence can be leveraged in the current ACOPF solvers to learn the iteration process and update parameters. Fig. \ref{rnn} shows the structure of the RNN algorithm. By using this method, the computation time may be reduced significantly and the final solutions can be guaranteed to be feasible and convergent.

\subsection{Choosing the Most Significant Inputs for ML Algorithms}

Although some current research is focused on learning the topology of ACOPF problems directly by training an ML model with existing data, researchers did not consider which inputs are necessary and which can be ignored. Some papers considered all variables in the ACOPF problem, and others used select variables based on their experience. However, no work systematically studied the significance of each variable. As can be seen from \eqref{acopf}, there are correlations among the voltages and active/reactive power flows in the ACOPF problem that allow one to  be denoted by the other. Therefore, there is potential to discover the importance of each variable and select the most significant ones as the inputs of the ML model. By using these essential and reduced data as the inputs, the training process may be improved effectively, and the model may work better in the predicting process and generate more convergent solutions.

\subsection{Discovering the Potential of Other ML Algorithms}

As shown in Section \uppercase\expandafter{\romannumeral3}, current research mainly utilized basic ML algorithms such as CNN, DNN, and so on. However, other ML methods may be more suitable for some OPF problems. For instance, for load-shedding and electric market price-prediction problems, LSTMs may be more appropriate since they are useful in time-series predictions based on previous inputs.

\begin{figure}[H]
\centering
\centerline{\includegraphics[width=0.6\textwidth]{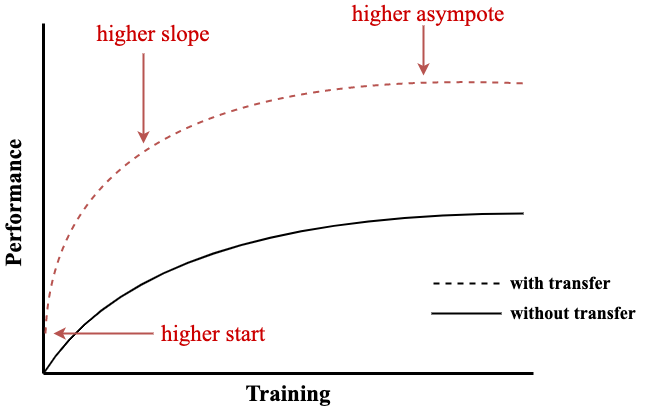}}
\caption{Three benefits of transfer learning.}
\label{tl}
\end{figure}

Another ML technique, transfer learning, aims to adapt a pretrained model on a new related target to improve learning performance and efficiency \cite{TLsurvey}. By combining transfer learning with traditional ML methods, three possible benefits may be achieved: 1) a higher initial performance; 2) a steeper slope (less time taken) during the learning process; and 3) a better final performance of the model trained \cite{HandbookOR}. A conceptual illustration of these three benefits is shown in Fig. \ref{tl}.

\begin{figure}[H]
\centering
\centerline{\includegraphics[width=0.6\textwidth]{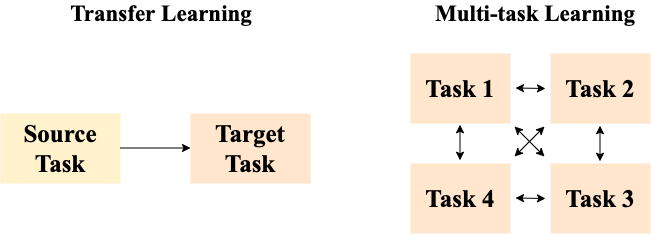}}
\caption{Distinction between transfer learning and multi-task learning.}
\label{dist}
\end{figure}

By taking advantage of these benefits in transfer learning, future research on ACOPF with ML can significantly improve the performance of traditional ML algorithms, or replace some learning methods with transfer learning, such as the MTL used in \cite{pgsim}. In MTL, the learning agent receives multiple tasks at one time, rather than be allocated with individual source and target tasks. That is, the information can be shared freely among all tasks. But in transfer learning, the agent has no information of a target task when it learns a source task since the information only flies in one direction (the distinction between transfer learning and MTL is shown in Fig. \ref{dist}). Obviously, multi-task learning is closely related to transfer learning, and it is possible to substitute multi-task learning with transfer learning, but reverse is impossible \cite{TLsurvey}. 

\subsection{Choosing the Most Appropriate ML Model }

As summarized in Section \uppercase\expandafter{\romannumeral3}, current research utilized popular ML algorithms for ACOPF problems in different situations, but a thorough comparison of these algorithms to determine which model fits the problem better is lacking. For example, the authors in \cite{LeveragingPG} evaluated the performances of fully connected neural network (FCNN), CNN and GNN models for two fundamental approaches to OPF with ML: regression (predicting optimal set points of generator) and classification (predicting the active set of constraints). The comparison of results showed that the FCNN models, which took less training time, are most appropriate for generator set-points and electricity prices problems. 

Inspired by \cite{LeveragingPG}, future research on ACOPF with ML algorithms should systematically compare different ML algorithms under different scenarios to provide clear guidance on the most appropriate model and improvements needed under various conditions. 

\subsection{Improving the Efficiency of ML and OPF Computation}

To improve model accuracy, ML algorithms require large amounts of data, requiring hours or even days for the training process. Also, OPF/ACOPF problems need to solve the power flow equations iteratively to obtain more optimal solutions, resulting in longer calculation times for larger power systems. Therefore, utilizing high-performance hardware to improve ML and OPF calculation efficiency should also be considered in future research. 

GPU (graphics processing unit), whose highly parallel structure is very friendly to algorithms that process large blocks of data in parallel, is more efficient than CPU (central processing units). Currently, most ML research utilizes GPU to tackle with the large number of data and some research on OPF problems also investigate the advantage of GPU \cite{gpu,AcceleratedCA}. In future research, combining ML algorithms with GPU as a co-processor to solve OPF/ACOPF problems is an intelligent and efficient method. Indeed, this promising combination on OPF/ACOPF problems really needs further investigation and development.  


\subsection{Improving the Validation \& Evaluation of ML Algorithms}

\subsubsection{Validating the ML Models on Larger Power Systems}
From Section \uppercase\expandafter{\romannumeral3}, although recent research on ACOPF with ML algorithms utilized multiple methods, their efficiency was validated predominantly  on <name> an up to $300$-bus systems. In reality, the number of buses far exceeds $300$, so future research should validate trained ML models on larger power systems to better evaluate their performance and applicability.

\subsubsection{Improving the Evaluation of the ML Algorithms}
As summarized in Table \ref{indices}, current research on ACOPF problems with ML algorithms only considered part of the metrics to evaluate the developed models and prediction performance, which may lack fairness and thoroughness. Therefore,  in the future researchers should consider all the related indices of the model evaluation. By doing this, a clearer path to improve NN models may emerge as well as facilitate comparisons with peer researchers. Also, it maybe beneficial to find the most appropriate ML model for the ACOPF problems. 

\section{Conclusion}

More than $50$ years after the nonlinear and nonconvex ACOPF problems were proposed, there are still no effective and exact methods that can obtain optimal solutions for them. Existing approaches can be classified into two main categories: conventional and AI techniques. However, traditional methods cannot guarantee the optimal solutions, and they are computationally burdensome. 

With the development of ML, many researchers are using this  technique to solve the nonconvex ACOPF problems. This paper reviews the most recent research on ACOPF with ML algorithms. By using ML, the accuracy and computation efficiency of the proposed algorithms were improved when compared with traditional methods. However, it is still an open question which algorithm is more suitable to solve ACOPF problems and which method can obtain the exact optimal solutions with high efficiency. This paper also discusses future research directions which include: 1) solving the relaxed ACOPF problems with ML algorithms; 2) improving current ACOPF solvers with ML algorithms; 3) choosing the most significant inputs of ML algorithms; 4) discovering the potential of other ML algorithms; 5) improving the efficiency of ML and OPF computation; and 6) improving the validation and evaluation of ML algorithms. In summary, research on ACOPF problems are still an active area of investigation all over the world, and ML can open a new door to future research on ACOPF.

\appendix

 \bibliographystyle{elsarticle-num} 
 \bibliography{cas-refs}

\begin{thebibliography}{10}
\expandafter\ifx\csname url\endcsname\relax
  \def\url#1{\texttt{#1}}\fi
\expandafter\ifx\csname urlprefix\endcsname\relax\def\urlprefix{URL }\fi
\expandafter\ifx\csname href\endcsname\relax
  \def\href#1#2{#2} \def\path#1{#1}\fi

\bibitem{first}
J.~Carpentier, {Contribution á l’étude du dispatching économique}, Bulletin de la Société Française des Électriciens 3 (1962) 431--447.

\bibitem{Frank2016AnIT}
S.~M. Frank, S.~Rebennack, {An Introduction to Optimal Power Flow: Theory, Formulation, and Examples}, IIE Transactions 48 (2016) 1172 -- 1197.

\bibitem{Zhao2021ASB}
M.~Zhao, M.~Barati, {A Sensitive-Eigenvalues based Global Algorithm for Penalized SDP Relaxation of ACOPF}, 2021 IEEE Power \& Energy Society Innovative Smart Grid Technologies Conference (ISGT) (2021) 1--5.

\bibitem{r1}
B.~C. {Lesieutre}, I.~A. {Hiskens}, {Convexity of the Set of Feasible Injections and Revenue Adequacy in {FTR} Markets}, IEEE Transactions on Power Systems 20~(4) (2005) 1790--1798.

\bibitem{history}
M.~B. Cain, R.~P. O'Neill, A.~Castillo, et~al., {History of Optimal Power Flow and Formulations}, Federal Energy Regulatory Commission 1 (2012) 1--36.

\bibitem{r3}
S.~{Frank}, S.~{Rebennack}, {An Introduction to Optimal Power Flow: Theory, Formulation, and Examples}, IIE Transactions 48~(12) (2016) 1172--1197.

\bibitem{r4}
M.~{Cain}, R.~{O'Neil}, A.~{Castillo}, {History of Optimal Power Flow and Formulations}, Federal Energy Regulatory Commission (2012) 1--36.

\bibitem{Jabr2007ACQ}
R.~A. Jabr, {A Conic Quadratic Format for the Load Flow Equations of Meshed Networks}, IEEE Transactions on Power Systems 22 (2007) 2285--2286.

\bibitem{Lavaei2012ZeroDG}
J.~Lavaei, S.~H. Low, {Zero Duality Gap in Optimal Power Flow Problem}, IEEE Transactions on Power Systems 27 (2012) 92--107.

\bibitem{conv}
M.~Ma, L.~Fan, Z.~Miao, B.~Zeng, H.~Ghassempour, {A Sparse Convex {AC OPF} Solver and Convex Iteration Implementation based on 3-node Cycles}, Electric Power Systems Research 180 (2020).

\bibitem{r2}
W.~A. {Bukhsh}, A.~{Grothey}, K.~I.~M. {McKinnon}, P.~A. {Trodden}, {Local Solutions of the Optimal Power Flow Problem}, IEEE Transactions on Power Systems 28~(4) (2013) 4780--4788.

\bibitem{20}
A.~S. {Zamzam}, X.~{Fu}, N.~D. {Sidiropoulos}, {Data-Driven Learning-Based Optimization for Distribution System State Estimation}, IEEE Transactions on Power Systems 34~(6) (2019) 4796--4805.

\bibitem{Baker}
K.~Baker, {Learning Warm-Start Points For {AC} Optimal Power Flow}, 2019 IEEE 29th International Workshop on Machine Learning for Signal Processing (MLSP) (2019) 1--6.

\bibitem{24}
D.~{Bertsimas}, B.~{Stellato}, {Online Mixed-Integer Optimization in Milliseconds}, ArXiv (2019).

\bibitem{25}
X.~{Pan}, T.~{Zhao}, M.~{Chen}, {{DeepOPF}: Deep Neural Network for {DC} Optimal Power Flow}, 2019 IEEE International Conference on Communications, Control, and Computing Technologies for Smart Grids (SmartGridComm) (2019) 1--6.

\bibitem{22}
S.~{Misra}, L.~Roald, Y.~{Ng}, {Learning for Constrained Optimization: Identifying Optimal Active Constraint Sets}, arXiv: Optimization and Control (2018).

\bibitem{23}
D.~{Deka}, S.~{Misra}, {Learning for {DC-OPF}: Classifying Active Sets using Neural Nets}, 2019 IEEE Milan PowerTech (2019) 1--6.

\bibitem{Guha}
N.~Guha, Z.~Wang, M.~Wytock, A.~Majumdar, {Machine Learning for {AC} Optimal Power Flow}, 2019 Climate Change Workshop at ICML (2019).

\bibitem{6756976}
S.~H. Low, Convex relaxation of optimal power flow—part i: Formulations and equivalence, IEEE Transactions on Control of Network Systems 1~(1) (2014) 15--27.
\newblock \href {https://doi.org/10.1109/TCNS.2014.2309732} {\path{doi:10.1109/TCNS.2014.2309732}}.

\bibitem{5971792}
J.~Lavaei, S.~H. Low, Zero duality gap in optimal power flow problem, IEEE Transactions on Power Systems 27~(1) (2012) 92--107.
\newblock \href {https://doi.org/10.1109/TPWRS.2011.2160974} {\path{doi:10.1109/TPWRS.2011.2160974}}.

\bibitem{7271127}
C.~Coffrin, H.~L. Hijazi, P.~Van~Hentenryck, The qc relaxation: A theoretical and computational study on optimal power flow, IEEE Transactions on Power Systems 31~(4) (2016) 3008--3018.
\newblock \href {https://doi.org/10.1109/TPWRS.2015.2463111} {\path{doi:10.1109/TPWRS.2015.2463111}}.

\bibitem{6980142}
D.~K. Molzahn, I.~A. Hiskens, Sparsity-exploiting moment-based relaxations of the optimal power flow problem, IEEE Transactions on Power Systems 30~(6) (2015) 3168--3180.
\newblock \href {https://doi.org/10.1109/TPWRS.2014.2372478} {\path{doi:10.1109/TPWRS.2014.2372478}}.

\bibitem{onil}
A.~O. Richard~P., Castillo, M.~B. Cain, The iv formulation and linear approximations of the ac optimal power flow problem, FERC (2012).

\bibitem{7299702}
A.~Castillo, P.~Lipka, J.-P. Watson, S.~S. Oren, R.~P. O'Neill, A successive linear programming approach to solving the iv-acopf, IEEE Transactions on Power Systems 31~(4) (2016) 2752--2763.
\newblock \href {https://doi.org/10.1109/TPWRS.2015.2487042} {\path{doi:10.1109/TPWRS.2015.2487042}}.

\bibitem{ml}
T.~Mitchell, {Machine Learning}, New York: McGraw Hill (1997).

\bibitem{mll}
datascience@berkeley, the online Master of Information and Data Science from UC Berkeley.

\bibitem{1}
A.~{Radovic}, {Neutrino Identification with a Convolutional Neural Network in the {NOvA} Detectors}, nternational Conference on High Energy Physics (2016).

\bibitem{2}
E.~{Racah}, C.~{Beckham}, T.~{Maharaj}, S.~{Kahou}, Prabhat, C.~{Pal}, {Extreme Weather: A Large-scale Climate Dataset for Semi-supervised Detection, Localization, and Understanding of Extreme Weather Events}, in: NIPS, 2017.

\bibitem{3}
W.~Dong, J.~Liu, Z.~Xie, D.~Li, {Adaptive Neural Network-based Approximation to Accelerate Eulerian Fluid Simulation}, Proceedings of the International Conference for High Performance Computing, Networking, Storage and Analysis (2019).

\bibitem{che}
R.~M. Balabin, E.~I. Lomakina, {Neural Network Approach to Quantum-Chemistry Data: Accurate Prediction of Density Functional Theory Energies.}, The Journal of chemical physics 131~(7) (2009).

\bibitem{game}
D.~Silver, et~al, {Mastering the Game of Go with Deep Neural Networks and Tree Search}, Nature 529~(7587) (2016) 484--489.

\bibitem{sc}
N.~Sengupta, M.~Sahidullah, G.~Saha, {Lung Sound Classification using Cepstral-based Statistical Features}, Computers in biology and medicine 75~(1) (2016) 118--129.

\bibitem{cancer}
N.~Ganesan, K.~Venkatesh, M.~A. Rama, A.~Palani, {Application of Neural Networks in Diagnosing Cancer Disease using Demographic Data}, International Journal of Computer Applications 1~(26) (2010) 81--97.

\bibitem{4}
L.~Richard, R.~Bonifetto, S.~Carli, A.~Froio, A.~Foussat, R.~Zanino, {Artificial Neural Network {(ANN)} modeling of the pulsed heat load during {ITER CS} magnet operation}, Cryogenics 63 (2014) 231--240.

\bibitem{Zamzam2020LearningOS}
A.~S. Zamzam, K.~Baker, {Learning Optimal Solutions for Extremely Fast {AC} Optimal Power Flow}, 2020 IEEE International Conference on Communications, Control, and Computing Technologies for Smart Grids (SmartGridComm) (2020) 1--6.

\bibitem{pgsim}
W.~Dong, Z.~Xie, G.~Kestor, D.~Li, {{Smart-PGSim}: Using Neural Network to Accelerate {AC-OPF} Power Grid Simulation}, SC20: International Conference for High Performance Computing, Networking, Storage and Analysis (2020) 1--15.

\bibitem{pred}
T.~Navidi, S.~Bhooshan, A.~Garg, {Predicting Solutions to the {Optimal Power Flow} Problem}, 2016.

\bibitem{Baker2020ALQ}
K.~Baker, {A Learning-boosted Quasi-Newton Method for AC Optimal Power Flow}, arXiv: Optimization and Control (2020).

\bibitem{MIPS}
R.~D. Zimmerman, H.~Wang, {MATPOWER Interior Point Solver (MIPS) User's manual}, Available: https://matpower.org/docs/ MIPS-manual-1.4.pdf, doi: 10.5281/zenodo.4073324 (Verision 1.4, 2020).

\bibitem{random}
G.~James, D.~Witten, T.~Hastie, R.~Tibshirani, {An Introduction to Statistical Learning: with Applications in {R}, Springer Texts in Statistics}, New York: Springer Science \& Business Media (2013).

\bibitem{DeepLM}
L.~{Deng}, D.~{Yu}, {Deep Learning: Methods and Applications}, Found. Trends Signal Process. 7 (2014) 197--387.

\bibitem{Zhou2020ADM}
Y.~Zhou, B.~Zhang, C.~Xu, T.~Lan, D.~Ruisheng, D.~Shi, Z.~Wang, W.~Lee, {A Data-driven Method for Fast {AC} Optimal Power Flow Solutions via Deep Reinforcement Learning}, Journal of Modern Power Systems and Clean Energy 8 (2020) 1128--1139.

\bibitem{ProximalPO}
J.~Schulman, F.~Wolski, P.~Dhariwal, A.~Radford, O.~Klimov, Proximal policy optimization algorithms, ArXiv abs/1707.06347 (2017).

\bibitem{RealTimeOP}
J.~H. Woo, L.~Wu, J.-B. Park, J.~Roh, {Real-Time Optimal Power Flow Using Twin Delayed Deep Deterministic Policy Gradient Algorithm}, IEEE Access 8 (2020) 213611--213618.

\bibitem{Singh2021LearningTS}
M.~Singh, V.~Kekatos, G.~Giannakis, {Learning to Solve the AC-OPF using Sensitivity-Informed Deep Neural Networks}, ArXiv abs/2103.14779 (2021).

\bibitem{Physicsin}
M.~Raissi, P.~Perdikaris, G.~E. Karniadakis, {Physics-informed Neural Networks: A Deep Learning Framework for Solving Forward and Inverse Problems Involving Nonlinear Partial Differential Equations}, Journal of Computational Physics 378 (2019) 686--707.

\bibitem{dnn}
F.~Fioretto, T.~Mak, P.~Hentenryck, {Predicting {AC} Optimal Power Flows: Combining Deep Learning and Lagrangian Dual Methods}, in: AAAI, 2020.

\bibitem{DeepOPFAF}
X.~Pan, M.~Chen, T.~Zhao, S.~Low, Deepopf: A feasibility-optimized deep neural network approach for ac optimal power flow problems, ArXiv abs/2007.01002 (2020).

\bibitem{PhysicsInformedNN}
R.~Nellikkath, S.~Chatzivasileiadis, {Physics-Informed Neural Networks for Minimising Worst-Case Violations in DC Optimal Power Flow}, ArXiv abs/2107.00465 (2021).

\bibitem{pypower}
Pypower, https://pypi.org/project/PYPOWER/ (2018).

\bibitem{2019ShortTermRL}
W.~Kong, Z.~Dong, Y.~Jia, D.~J. Hill, Y.~Xu, Y.~Zhang, {Short-Term Residential Load Forecasting Based on LSTM Recurrent Neural Network}, IEEE Transactions on Smart Grid 10 (2019) 841--851.

\bibitem{bilstm}
C.~Cai, Y.~Tao, T.~Zhu, Z.~Deng, {Short-Term Load Forecasting Based on Deep Learning Bidirectional LSTM Neural Network}, Applied Sciences 11~(17) (2021) 8129.

\bibitem{cnnlstm}
S.~H. Rafi, N.~Masood, S.~R. Deeba, E.~Hossain, {A Short-Term Load Forecasting Method Using Integrated CNN and LSTM Network}, IEEE Access 9 (2021) 32436--32448.

\bibitem{Kim2019GraphCN}
C.~Kim, K.~Kim, P.~Balaprakash, M.~Anitescu, {Graph Convolutional Neural Networks for Optimal Load Shedding under Line Contingency}, 2019 IEEE Power \& Energy Society General Meeting (PESGM) (2019) 1--5.

\bibitem{BayesianSE}
K.~Mestav, J.~Luengo-Rozas, L.~Tong, {Bayesian State Estimation for Unobservable Distribution Systems via Deep Learning}, IEEE Transactions on Power Systems 34 (2019) 4910--4920.

\bibitem{2021ModelingTA}
A.~Kody, S.~C. Chevalier, S.~Chatzivasileiadis, D.~K. Molzahn, {Modeling the AC Power Flow Equations with Optimally Compact Neural Networks: Application to Unit Commitment}, ArXiv abs/2110.11269 (2021).

\bibitem{GraphNN}
S.~Liu, C.~Wu, H.~Zhu, {Graph Neural Networks for Learning Real-Time Prices in Electricity Market}, ArXiv abs/2106.10529 (2021).

\bibitem{2018UsingML}
L.~Duchesne, E.~Karangelos, L.~Wehenkel, {Using Machine Learning to Enable Probabilistic Reliability Assessment in Operation Planning}, 2018 Power Systems Computation Conference (PSCC) (2018) 1--8.

\bibitem{2020RecentDI}
L.~Duchesne, E.~Karangelos, L.~Wehenkel, {Recent Developments in Machine Learning for Energy Systems Reliability Management}, Proceedings of the IEEE 108 (2020) 1656--1676.

\bibitem{InverseOP}
P.~Donti, I.~Azevedo, J.~Kolter, {Inverse Optimal Power Flow: Assessing the Vulnerability of Power Grid Data}, 2018.

\bibitem{GaussianPL}
P.~Pareek, H.~Nguyen, {Gaussian Process Learning-Based Probabilistic Optimal Power Flow}, IEEE Transactions on Power Systems 36 (2021) 541--544.

\bibitem{Zhao2021LoworderMR}
M.~Zhao, M.~Barati, {Low-order Moment Relaxation of ACOPF via Algorithmic Successive Linear Programming}, 2021 IEEE Texas Power and Energy Conference (TPEC) (2021) 1--6.

\bibitem{Shen2018APM}
X.~Shen, J.~E. Mitchell, {A Penalty Method for Rank Minimization Problems in Symmetric Matrices}, Computational Optimization and Applications 71 (2018) 353--380.

\bibitem{TLsurvey}
F.~Zhuang, Z.~Qi, K.~Duan, D.~Xi, Y.~Zhu, H.~Zhu, H.~Xiong, Q.~He, {A Comprehensive Survey on Transfer Learning}, Proceedings of the IEEE 109 (2021) 43--76.

\bibitem{HandbookOR}
E.~S. Olivas, J.~D.~M. Guerrero, M.~M. Sober, J.~Benedito, A.~L{\'o}pez, {Handbook Of Research On Machine Learning Applications and Trends: Algorithms, Methods and Techniques (2 Volumes)}, 2009.

\bibitem{LeveragingPG}
T.~Falconer, L.~Mones, {Leveraging Power Grid Topology in Machine Learning Assisted Optimal Power Flow}, ArXiv abs/2110.00306 (2021).

\bibitem{gpu}
Z.~Chen, J.~Liu, X.~Liu, {GPU Accelerated Power Flow Calculation of Integrated Electricity and Heat System with Component-oriented Modeling of District Heating Network}, Applied Energy 305 (2022).

\bibitem{AcceleratedCA}
Y.~Kim, K.~Kim, {Accelerated Computation and Tracking of AC Optimal Power Flow Solutions using GPUs}, ArXiv 2110.06879 (2021).

\end{thebibliography}





\end{document}